\documentclass{amsart}

\usepackage{amssymb}	
\usepackage{tikz}
\usepackage{tikz-cd}

\usepackage{amsmath}
 
\usepackage{mathrsfs}
\usepackage{amsthm}
\usepackage{verbatim}
\usepackage{subeqnarray}
\usepackage{graphicx}
\usepackage{cancel}
 
\usepackage{color}
 
\usepackage{float}
 
\usepackage{bm}
\usepackage{hyperref}
 
\usepackage{amsfonts}
 
\usepackage{amscd}
\usepackage{cases}
\usepackage{xcolor}
\usepackage{eucal}

\definecolor{red}{rgb}{0.8,0,0}
\definecolor{darkorange}{rgb}{1,0.4,0}
\definecolor{lightorange}{rgb}{1,0.6, 0}
\definecolor{yellow}{rgb}{1,0.8, 0}

\usepackage{chngpage}

\usepackage{lipsum}

\usepackage{booktabs}

\usepackage{geometry}
\begin{document}
\title{Oberwolfach report : Discretization of Hilbert complexes}
 \author{Kaibo Hu}
 \address{Mathematical Institute,
University of Oxford, Andrew Wiles Building, Radcliffe Observatory Quarter, 
Oxford, OX2 6GG, UK}
\email{Kaibo.Hu@maths.ox.ac.uk}
\date{}

\begin{abstract}
The report is based on an extended abstract for the MFO workshop ``Hilbert Complexes: Analysis, Applications, and Discretizations'', held at Oberwolfach during 19-25 June 2022. The aim is to provide an overview of some aspects of discretization of Hilbert complexes with an emphasis on conforming finite elements. 
\end{abstract}

\maketitle
 
Discretization of Hilbert complexes plays a central role in finite element exterior calculus \cite{arnold2018finite,Arnold.D;Falk.R;Winther.R.2006a,Arnold.D;Falk.R;Winther.R.2010a}. In this report, we review some recent development in this direction, emphasizing conforming finite elements.

 \medskip
{\noindent \bf Overview.} The overall aim of discretizing Hilbert complexes is to construct finite-dimensional spaces that fit in a complex. These spaces should have certain continuity and inherit algebraic and differential structures on the continuous level, e.g., exactness. The construction of spaces with high regularity is naturally connected with spline theory. In contrast, allowing low regularity often leads to spaces involving Dirac measures, which we will refer to as {\it distributional finite elements} \cite{braess2008equilibrated,christiansen2011linearization,licht2017complexes}. The construction of simplicial splines is subtle due to intrinsic supersmoothness \cite{floater2020characterization,sorokina2010intrinsic} (piecewise smooth functions may have automatic higher order continuity at hinges, e.g., vertices and edges of a mesh). The construction of distributional finite elements is usually more straightforward, but numerical schemes call for attention. The general idea is to choose proper combinations of spaces such that Dirac measures are only evaluated on continuous functions or forms \cite{gopalakrishnan2020mass,pechstein2011tangential}. The situation can be more subtle for nonlinear problems since the product of Dirac measures is not defined in general (c.f. \cite{berchenko2021finite,christiansen2013exact,gopalakrishnan2022analysis}).

This report will mainly focus on the construction of conforming finite elements. On cubical meshes, the tensor product provides a useful algebraic tool to derive elements and complexes in any dimension from results in 1D \cite{arnold2015finite,cubicBGG,bonizzoni2021h,buffa2011isogeometric}. 
 On simplicies, there are two general strategies to handle supersmoothness: incorporating supersmoothness in the definition of finite element spaces (e.g., the Argyris $C^{1}$ element with second order vertex derivative degrees of freedom) or subdividing a simplex and seeking piecewise polynomials on the refined mesh (e.g., the Clough-Tocher $C^{1}$ element). In general, several possible versions exist and finer subdivisions usually require less supersmoothness, c.f. \cite{christiansen2018generalized,fu2020exact,guzman2022exact,neilan2015discrete}.
 
 \medskip
 {\noindent\it \bf De-Rham complex.}
The development of discrete de-Rham complexes demonstrates deep connections between topology and numerical analysis. The Raviart-Thomas \cite{Raviart.P;Thomas.J.1977a}, N\'ed\'elec \cite{Nedelec.J.1980a,Nedelec.J.1986a} and Brezzi-Douglas-Marini \cite{brezzi1985two} finite elements achieved success in computational electromagnetism and other areas. It was realized later that these finite elements fit in a discrete de-Rham complex, and the lowest order version coincides with the Whitney forms in geometric integration theory \cite{whitney2012geometric}. This observation inspired the development of discrete differential forms \cite{bossavit1988whitney,hiptmair1999canonical,hiptmair2001higher} and the finite element exterior calculus \cite{arnold2018finite,Arnold.D;Falk.R;Winther.R.2006a,Arnold.D;Falk.R;Winther.R.2010a}. The finite element periodic table \cite{arnold2014periodic} summarized several families of de-Rham finite elements on simplicial and cubical meshes. Poincar\'e (Koszul) operators, $\mathfrak{p}^{j}, \, j=1, 2, \cdots, n,$ satisfying $d^{k-1}\mathfrak{p}^{k}+\mathfrak{p}^{k+1}d^{k}=I$, provide explicit forms of potential and thus lead to polynomial exact sequences on cells \cite{Arnold.D;Falk.R;Winther.R.2006a,hiptmair1999canonical}. Geometric decomposition \cite{arnold2009geometric,chen2021geometric} provides insights for constructing degrees of freedom.

 \medskip
 {\noindent\it \bf Smoother de-Rham (Stokes) complexes.} Stokes problems in fluid mechanics raise the question of constructing a finite element velocity space $V_{h}\subset [H^{1}]^{3}$ and a pressure space $Q_{h}\subset L^{2}$ such that $\operatorname{div} V_{h}=Q_{h}$ and the inf-sup condition holds. The Stokes problem motivates discretizing de-Rham complexes that are smoother than the standard $H\Lambda$ version ($u$ in $L^{2}$ with $du$ in $L^{2}$). On the continuous level, there are several possible combinations of Sobolev spaces that fit in a de-Rham complex \cite{costabel2010bogovskiui}. Among them, the following has a direct application in the 2D Stokes problem \cite{falk2013stokes}:
 \begin{equation*}
 \begin{tikzcd}
 0\arrow{r}{}& H^{2}\arrow{r}{\operatorname{curl}}&{[H^{1}]}^{2}\arrow{r}{\operatorname{div}} &L^{2}\arrow{r}{} & 0.
 \end{tikzcd}
 \end{equation*}
The above sequence shows connections between scalar high order problems ($H^{2}$) and the Stokes problem ($\operatorname{div}: {[H^{1}]}^{2} \to L^{2}$). To a large extent, these two topics were developed independently over a long time (see, e.g., \cite{arnold1992quadratic,lai2007spline,zhang2005new,zhang2008p1}) until smoother de-Rham complexes built a bridge. On the one hand, one may differentiate scalar splines and obtain a Stokes pair. On the other hand, investigating the pre-image of the divergence operator helps to clarify the inf-sup stability of Stokes pairs (e.g., the Scott-Vogelius elements \cite{guzman2017scott,scott1985norm}).

Most existing conforming elements either use subdivision (e.g., \cite{arnold1992quadratic,fu2020exact,guzman2022exact,zhang2008p1}) or/and higher order derivatives as degrees of freedom (e.g., \cite{falk2013stokes,neilan2015discrete,neilan2020stokes}), which are two general strategies for handling supersmoothness.
Recent development in 2D and 3D can be found in, e.g., \cite{arf2021structure,christiansen2018generalized,guzman2020exact,guzman2018inf,hu2022family}. In the direction of reducing the number of degrees of freedom, one is interested in a minimal Stokes pair, where the degrees of freedom for the velocity space involve vertex evaluation (for approximation) and normal component on each face (for showing surjectivity of $\operatorname{div}$ with Stokes' theorem). 

 \centerline{
 \includegraphics[width=0.25\textwidth]{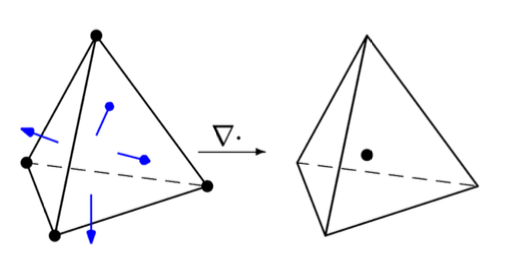} 
 }
{\noindent}Bernardi and Raugel \cite{bernardi1985conforming} constructed elements with the above degrees of freedom by enriching Lagrange elements with face bubbles. The resulting elements do not fit in a complex. Guzm\'an and Neilan \cite{guzman2018inf} modified the bubbles using modes on the Alfeld split of a tetrahedron to render the complex property (see also \cite{hu2022family} for the entire complex). A canonical construction on various subdivisions can be found in \cite{christiansen2018generalized}.

 \medskip
 {\noindent\it \bf BGG complexes.} 
 BGG complexes, as well as the BGG diagrams that lead to them, encode structures of many PDE problems \cite{arnold2021complexes,vcap2022bgg,pauly2020elasticity,pauly2020divdiv}. Since the construction of the Arnold-Winther element \cite{arnold2002mixed}, which was the first conforming triangular finite element with polynomial shape functions, several discretizations for linear elasticity were developed (e.g., \cite{hu2015family,hu2016finite}). Recently, there has been a surge of finite elements for Hilbert complexes, especially the Hessian, elasticity and $\operatorname{div}\operatorname{div}$ complexes in 2D and 3D \cite{arf2021structure,chen2020discrete,chen2020finite,chen2021finite,chen2021finite2,chen2021geometric,christiansen2020discrete,christiansen2022finite,hu2021conforming,hu2021conforming2,hu2022new,sander2021conforming}. 

The approach in \cite{arnold2021complexes} for deriving BGG complexes on the continuous level is to collect several copies of de-Rham complexes and eliminate some components. Therefore a natural approach to deriving discrete BGG complexes is to mimic the construction and fit finite element spaces in diagrams. A simple example is demonstrated in the following diagram, where one connects two 1D de-Rham complexes 
\begin{center}
\includegraphics[width=1.5in]{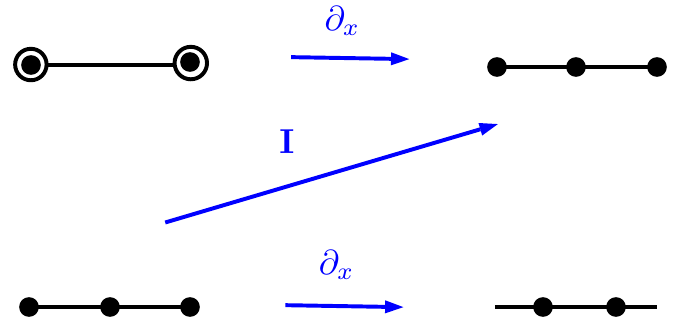} 
\end{center}
to derive a finite element complex with a second order differential operator:
\begin{center}
\includegraphics[width=1.5in]{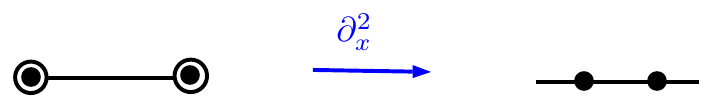}.
\end{center}
This example also demonstrates that the input discrete de-Rham complexes should have different regularity in this approach. 
Such finite element diagram chasing was first used in \cite{Arnold2006a} as a re-interpretation of the Arnold-Winther elasticity element. Recently, this approach has been extended to derive more complexes with applications in discretizing elasticity and curvature \cite{christiansen2020discrete,christiansen2018nodal,christiansen2022finite}. 

One may also directly discretize the BGG complexes without referring to the BGG diagrams. Hu and collaborators \cite{hu2021conforming,hu2021conforming2,hu2022new,hu2015family} constructed several conforming finite element complexes in 2D and 3D based on explicit characterizations of bubble functions and an investigation of Lagrange type bases. Chen and Huang \cite{chen2020discrete,chen2021geometric,chen2022finite-2D,chen2021finite,chen2021finite2,chen2020finite} further developed geometric decomposition and polynomial BGG complexes and obtained various conforming finite element complexes. There has been progress towards a systematic discretization of BGG complexes in any dimension \cite{cubicBGG,chen2021finite2} and a wide variety of continuity \cite{chen2022finite-2D,hu2021construction}.

 \medskip
 {\noindent\it\bf Summary and overlook.} This report aims to review some progress on discretizing Hilbert complexes. The emphasis is mainly on conforming finite elements on simplices, while other important topics are not covered, e.g., polyhedral elements, virtual elements, isogeometric analysis, nonconforming elements and applications. 
 
 Canonical finite element de-Rham complexes have been implemented in several packages, e.g., FEniCS \cite{fenics}, Firedrake \cite{firedrake}, NGSolve \cite{ngsolve}. Nevertheless, to the best of the author's knowledge, not many smoother finite element de-Rham (Stokes) complexes or BGG complexes have   been included. An exception is the Regge element, which has been implemented in several finite element packages \cite{fenics,firedrake,ngsolve}. Regge calculus  \cite{regge1961general} was proposed as a scheme for quantum and numerical gravity, and later interpreted as a finite element \cite{christiansen2011linearization,li2018regge}.  The Regge element provides another demonstration of the interactions between discrete theories and finite elements.
  
 Computational issues, e.g., well-conditioned bases and implementation, for smoother de-Rham and the BGG complexes require further investigation. Ideas and algorithms, e.g., Bernstein-Bezi\'er techniques, for scalar splines may be generalized to problems of vectors and tensors \cite{alfeld2016linear,lai2007spline,sorokina2018bernstein}. The algebraic structures in complexes may also provide a new perspective for theoretical questions of splines, e.g., dimension of spline spaces \cite{schenck2016multivariate}. General constructions of distributional and nonconforming elements for the BGG complexes call for further investigation. These elements enjoy simple degrees of freedom, and may thus provide a bridge for discretization of PDEs and discrete structures, including graph theory \cite{lim2020hodge}, discrete mechanics \cite{hauret2007diamond}, discrete differential geometry \cite{berchenko2021finite,christiansen2011linearization,christiansen2022finite,gopalakrishnan2022analysis,li2018regge} and gauge theory \cite{christiansen2012simplicial}. Finite elements may provide a new perspective for these areas by supplying local shape functions \cite{christiansen2011linearization} and inspire new schemes. Recent progress on Hilbert complexes paves a way to tackle the Einstein equations with applications in numerical relativity \cite{beig2020linearised,li2018regge,quenneville2015new}. 
 
\section*{Acknowledgement}
 
 The author would like to thank  Ana M. Alonso Rodriguez, Douglas N. Arnold, Dirk Pauly and Francesca Rapetti, for organizing the event, and other participants for stimulating discussions. The author thanks MFO for the excellent research facilities and Qian Zhang for proofreading the report. 
 
\bibliographystyle{siam}      
\bibliography{reference}{}   

\end{document}